\documentclass[reqno,tbtags,intlimits,a4paper,oneside,11pt]{amsart}
\usepackage{amsfonts}
\usepackage{}
\usepackage{amsfonts,amsmath,amssymb}
\newtheorem{Theorem}{Theorem}

\def\beq#1#2\eeq{%
        \begin{equation}%
        \label{#1}%
            #2%
        \end{equation}%
    }

\usepackage{color}
\usepackage{graphicx}
\usepackage{epstopdf}

\title[On shadow Mordell triples and Conway topograph]{Conway's light on the shadow of Mordell}

\author{A.P. Veselov}
\address{Department of Mathematical Sciences,
Loughborough University, Loughborough LE11 3TU, UK}
\email{A.P.Veselov@lboro.ac.uk}

\begin{document}

\maketitle

\begin{abstract}
Recently Valentin Ovsienko introduced a ``shadow" version of the celebrated Markov triples as the solutions of certain version of Markov equation over dual numbers.
We will discuss similar question for the Mordell Diophantine equation $$
X^2+Y^2+Z^2=2XYZ+1.
$$ We will see that the shadows of the special solution $(1,1,1)$ can be described using the original Conway topograph of the values of binary quadratic forms. The shadows of other Mordell triples will be written explicitly in terms of the corresponding solutions of Pell's equation. Their growth along the paths on the Conway topograph is described in terms of the Lyapunov function of the Euclid tree.
\end{abstract}


\section*{Introduction}

The celebrated Markov triples are positive integer solutions of the Markov equation
$$x^2 + y^2 + z^2=3xyz.$$ Markov \cite{Markov} showed that all such solutions can be found recursively from $(1,1,1)$ using a natural action of the modular group $PSL_2(\mathbb Z)$ on the solutions of this equation generated by the cyclic permutation of variables and Vieta involution $(x,y,z)\to (x,y, 3xy-z)$.

Markov triples were introduced by A.A. Markov in 1879 in relation with the Diophantine properties of the binary quadratic forms, which can be interpreted as the description of the ``most irrational" numbers. Since then there were discovered many surprising connections of Markov triples with hyperbolic geometry, combinatorics of words and, more recently, with the theory of Frobenius manifolds, algebraic geometry and cluster mutations (see \cite{Aigner, BV} and references therein). This makes Markov theme one of the most fascinating stories in mathematics, which because of still unproven uniqueness conjecture \cite{Aigner} is far from being finished.

Very recently Ovsienko \cite{O1}, motivated by his attempts to understand possible super analogues of cluster algebras \cite{O2, OS, OT},  introduced the so-called ``shadow" version of Markov triples by considering the solutions of a certain version of Markov equation over {\it dual numbers} $X=a+b\varepsilon, \varepsilon^2=0.$ 

The dual numbers were introduced in 1873 by Clifford and were used by Study to describe the relative position of two skew lines in space. They are the simplest case of the Grassmann numbers, which are used to describe the fermionic fields in modern quantum theory.

In fact, Ovsienko introduced the shadow versions of other sequences of integers including Fibonacci numbers and (with Tabachnikov) Somos-4 sequences \cite{O1, OT}.
Recently Hone \cite{Hone} extended this study to a more general Somos-4 recurrence relation 
over dual numbers.

As always with a new development, a valid question is whether the shadow dynamics is interesting enough to deserve a detailed study.
The aim of this note is to add few more arguments to Ovsienko and Hone in favour of the positive answer to this question. 

Namely, we will study the shadow Mordell triples, satisfying the {\it Mordell equation}
$$
X^2+Y^2+Z^2=2XYZ+1
$$
with $X,Y,Z$ being dual integers. As it was shown by Mordell \cite{Mordell}, this equation can be viewed as a ``solvable" modification of Markov equation. This fact was used by Zagier \cite{Zagier} to study the growth of the Markov numbers.

Note that, like in the case of Markov equation, we have natural action of the modular group $PSL_2(\mathbb Z)$ on the solutions of the Mordell equation. Remarkably the shadow orbits of the solution $(1,1,1)$, which is stationary for modular dynamics, are none other than the values of binary quadratic forms and thus can be visualised using the original Conway topograph \cite{Conway}. 

The shadows of other Mordell triples can be written explicitly in terms of the solutions of Pell's equations, following Mordell \cite{Mordell}. We will discuss also the shadow growth along the paths on the Conway topograph, which can be described in terms of the Lyapunov function of the Euclid tree \cite{SV1}.

\section*{Conway topograph and the modular group}

The values of binary quadratic forms on the integer lattice $ \mathbb{Z}^2$ is one of the most classical objects of study in number theory, going back to Fermat with later contributions from Euler, Gauss and Jacobi.
A famous question is: which integers $n$ can be represented as the sum of two squares:
$$
n=x^2+y^2, \, x,y \in \mathbb Z,
$$
and in how many ways. The answer is given by the Jacobi formula for the number $N(n)$ of such representations:
$$
N(n) = 4(N_1(n)
- N_3(n)),
$$
where $N_1(n)$ and $N_3(n)$ are the numbers of the
divisors of $n$ with the residues $1$ and $3$ modulo 4
respectively (see e.g. \cite{LeVeque}).

In 1990s this very classical story had a very interesting twist due to the intervention of John H. Conway.
In his book ``The Sensual (Quadratic) Form" Conway \cite{Conway} described the following ``topographic" way to ``visualise" these values.

Conway introduced the notions of the {\it lax} vector as a pair $(\pm v), v \in \mathbb Z^2$, and of the \emph{superbase} of the integer lattice $\mathbb Z^2$ as a triple of lax vectors $(\pm e_1, \pm e_2, \pm e_3)$ such that $(e_1, e_2)$ is a basis of the lattice and
		\begin{equation*}
		e_1 + e_2 + e_3 =  0.
		\end{equation*}
		
		Every basis $(e_1,e_2)$ can be included in exactly two superbases, namely  $\pm(e_1,e_2, -e_1-e_2)$ and $\pm(e_1,-e_2, -e_1+e_2)$, so that the corresponding set of the superbases can be described using the binary tree embedded in the plane.
The connected components of the complement to the tree (domains) are labelled by the primitive (that is, having coprime coordinates) lax vectors, the edges by the lax bases, while the superbases correspond to the vertices of the tree (see left hand side of Fig. 1, where we showed only one representative of the lax vectors).
The projective version of the Conway superbase topograph is known also as {\it Farey tree} since it is related to Farey ``addition" $\frac{a}{b}*\frac{c}{d}=\frac{a+c}{b+d}$ (see Fig. 1).

	
%
		
\begin{figure}[h]
\begin{center}
\includegraphics[height=48mm]{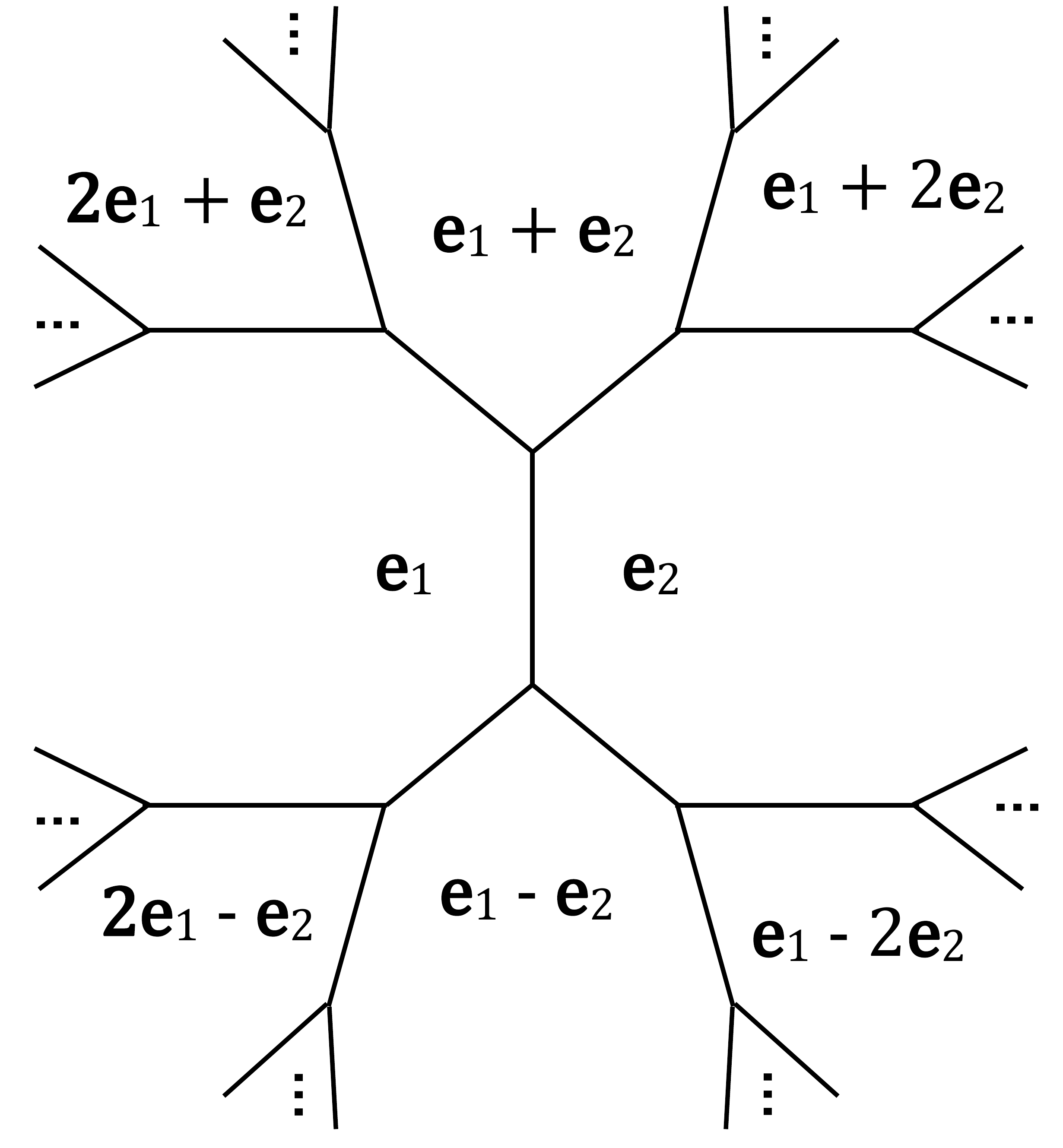}  \hspace{8pt}  \includegraphics[height=48mm]{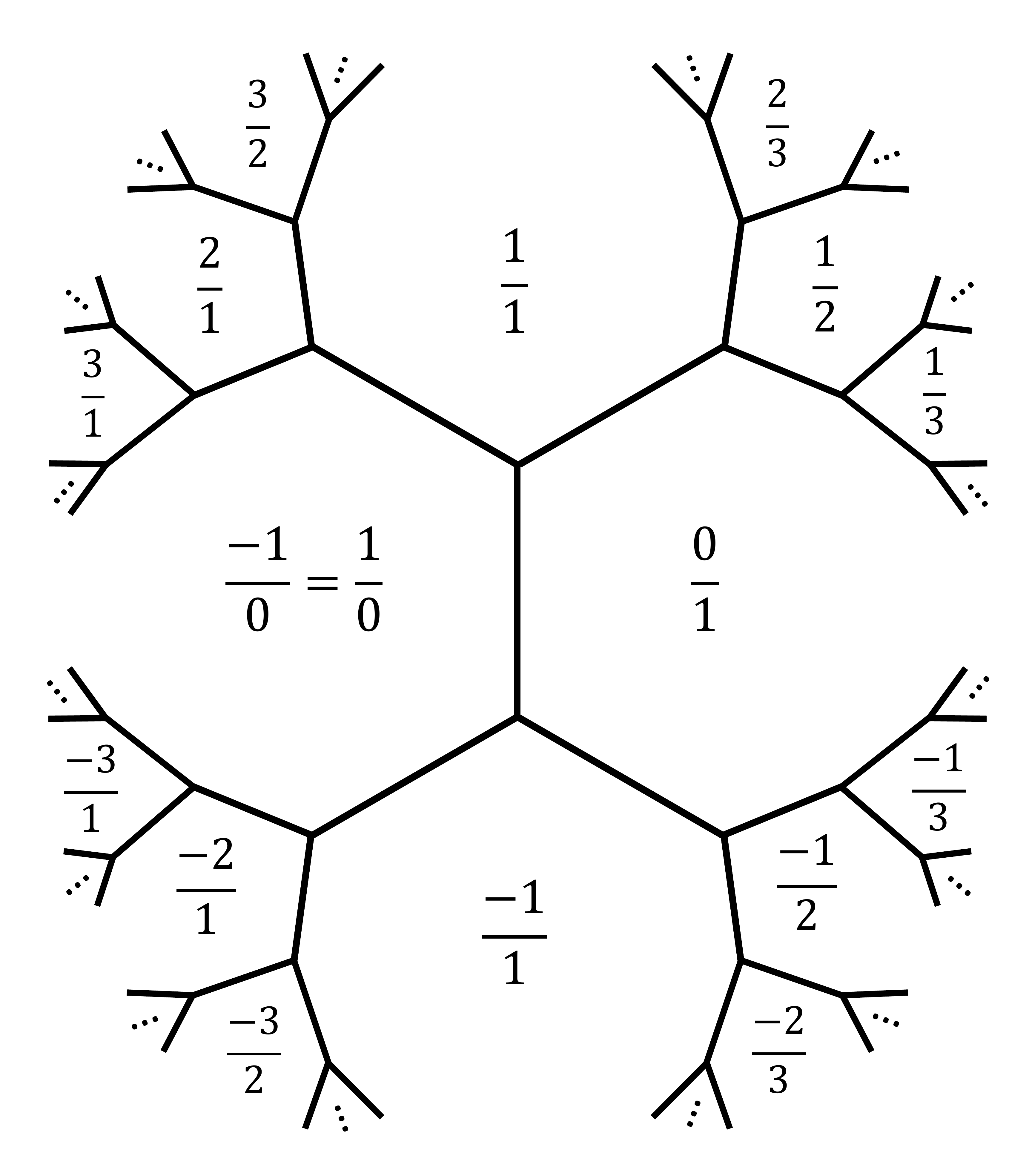}
\caption{\small The superbase topograph and the Farey tree.}
\end{center}
\end{figure}

%
%
	
Let
$Q(x, y) = ax^2 + hxy + by^2$
be a binary quadratic form (in Conway's notations), where $x,y \in \mathbb Z,$ but we can allow the coefficients $a,b,h$ to be not necessarily integer.

By taking the values of the form $Q$ on the vectors of the superbase, we get what Conway called the {\it topograph} of $Q$ containing the values of $Q$ on all primitive lattice vectors.
In particular, if $\mathbf{e}_1=(1,0), \mathbf{e}_2=(0,1), \mathbf{e}_3=-(1,1)$ we have the values 
$
Q(\mathbf{e}_1)=a, \,\, Q(\mathbf{e}_2)=b, \,\, Q(\mathbf{e}_3)=c:=a+b+h.
$

The key idea of Conway\footnote{Conway clearly believed that this is an important idea claiming in the Introduction to his lectures \cite{Conway} that ``the ``topograph" of the First Lecture is new."} is that one can construct the topograph of $Q$ {\it recursively} starting from this triple and using the following property of any quadratic form, which Conway called the {\it arithmetic progression rule} (known also in geometry as the {\it parallelogram law}):
\begin{equation}
\label{apr}
Q(\mathbf{u}+\mathbf{v})+Q(\mathbf{u}-\mathbf{v})=2(Q(\mathbf{u})+Q(\mathbf{v})), \quad \mathbf{u},\mathbf{v} \in \mathbb R^2.
\end{equation}
%
In fact, the quadratic forms (in any dimension) can be characterised as the continuous functions $Q$ satisfying relation (\ref{apr}).

In particular, for 
$
Q(x, y) = x^2 + xy + y^2, \, (x,y) \in \mathbb{Z}^2
$
describing the square lengths of vectors in regular hexagonal lattice, we have the Conway topograph shown on Fig. 2.

\begin{figure}[h]
	\centering
	\includegraphics[scale=0.2]{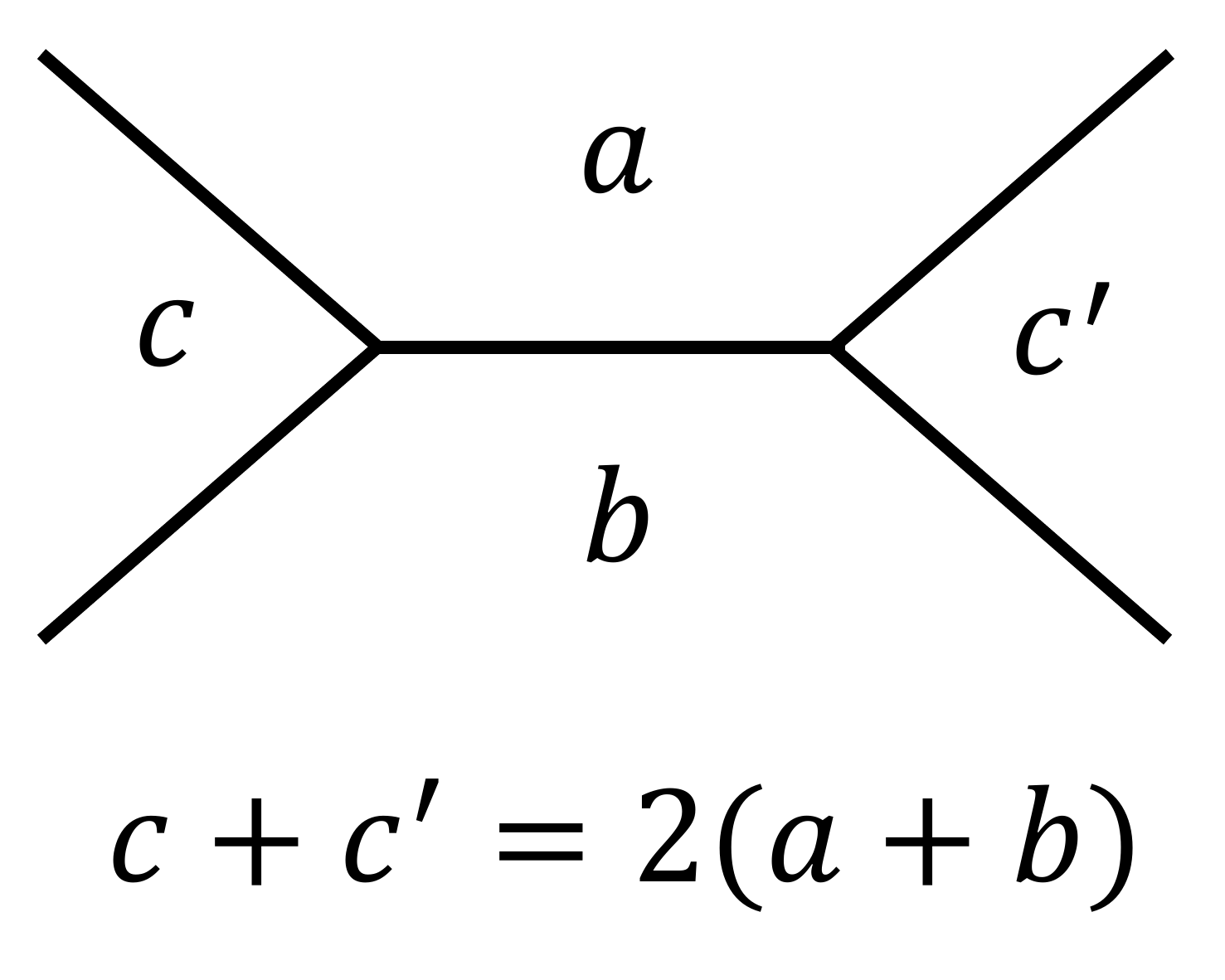} \quad \includegraphics[scale=0.5]{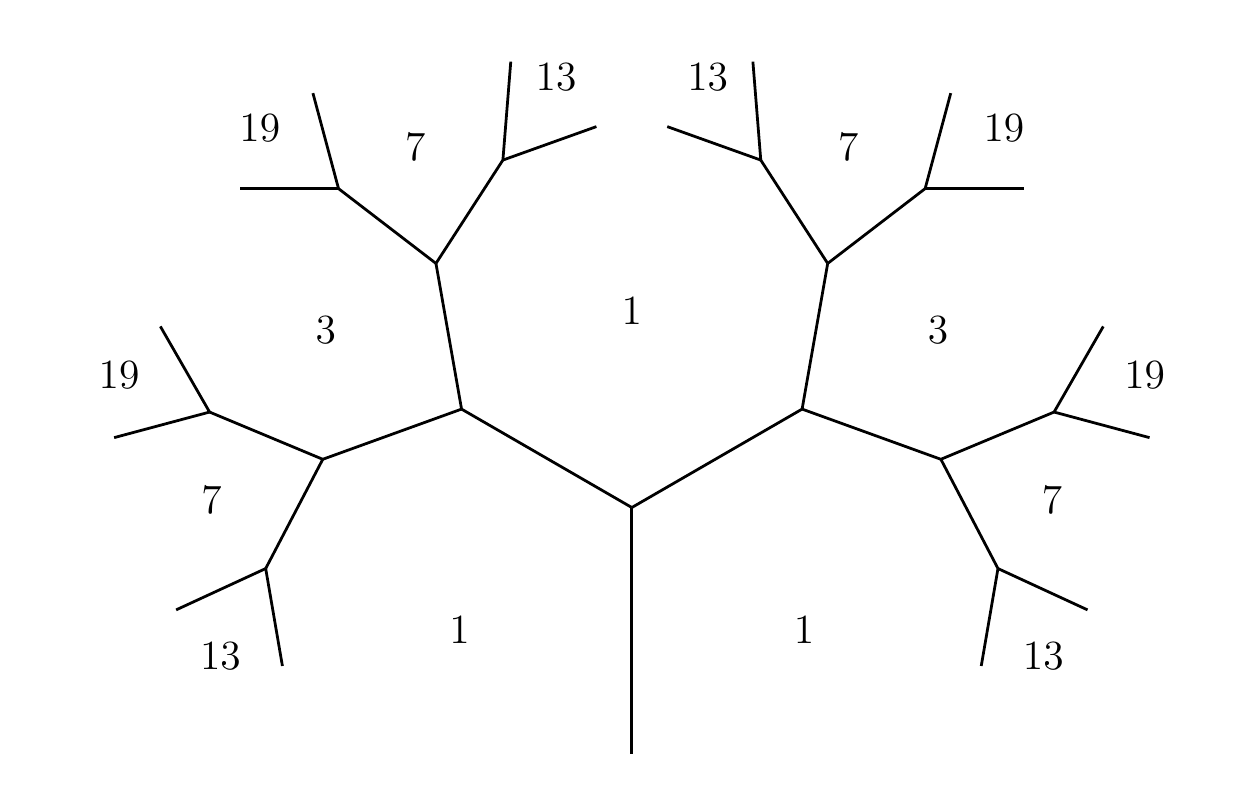} 
	\caption{Arithmetic progression rule and Conway topograph of quadratic form $Q=x^2+xy+y^2$}
	\end{figure}

As Conway nicely explained, the modular group $PSL_2(\mathbb Z)$ can be viewed as the natural symmetry group of his topograph. 
Indeed, it is well-known that $PSL(2, \mathbb Z)=\mathbb Z_2*\mathbb Z_3$ is generated by
$$
U=\begin{pmatrix}
  0 & -1 \\
 1 & 0 \\
\end{pmatrix},
V=\begin{pmatrix}
  0 & -1 \\
 1 & 1 \\
\end{pmatrix},
$$
which act on the topograph combinatorially by rotation by $\pi$ about edge centre and by $2\pi/3$ about vertex respectively. 
Thus, within Klein's Erlangen programme, Conway topograph can be viewed as a discrete version of the hyperbolic plane with its isometry group $PSL_2(\mathbb R)$ replaced by $PSL_2(\mathbb Z)$.

This means that Conway topograph can be used to describe any $PSL_2(\mathbb Z)$-dynamics. In particular, we have a natural action of $PSL_2(\mathbb Z)$ on the solutions of Markov equation
$x^2 + y^2 + z^2=3xyz$, generated by the Vieta involution 
 and cyclic permutation of the variables. As Markov showed in \cite{Markov}, all integer solutions of Markov equation form one orbit of $PSL_2(\mathbb Z)$ acting on $(1,1,1)$.
 A part of the corresponding Conway topograph of the Markov triples is shown on Fig. 3 with the local Vieta rule generating them.
 
 \begin{figure}[h]
\begin{center}
 \includegraphics[width=40mm]{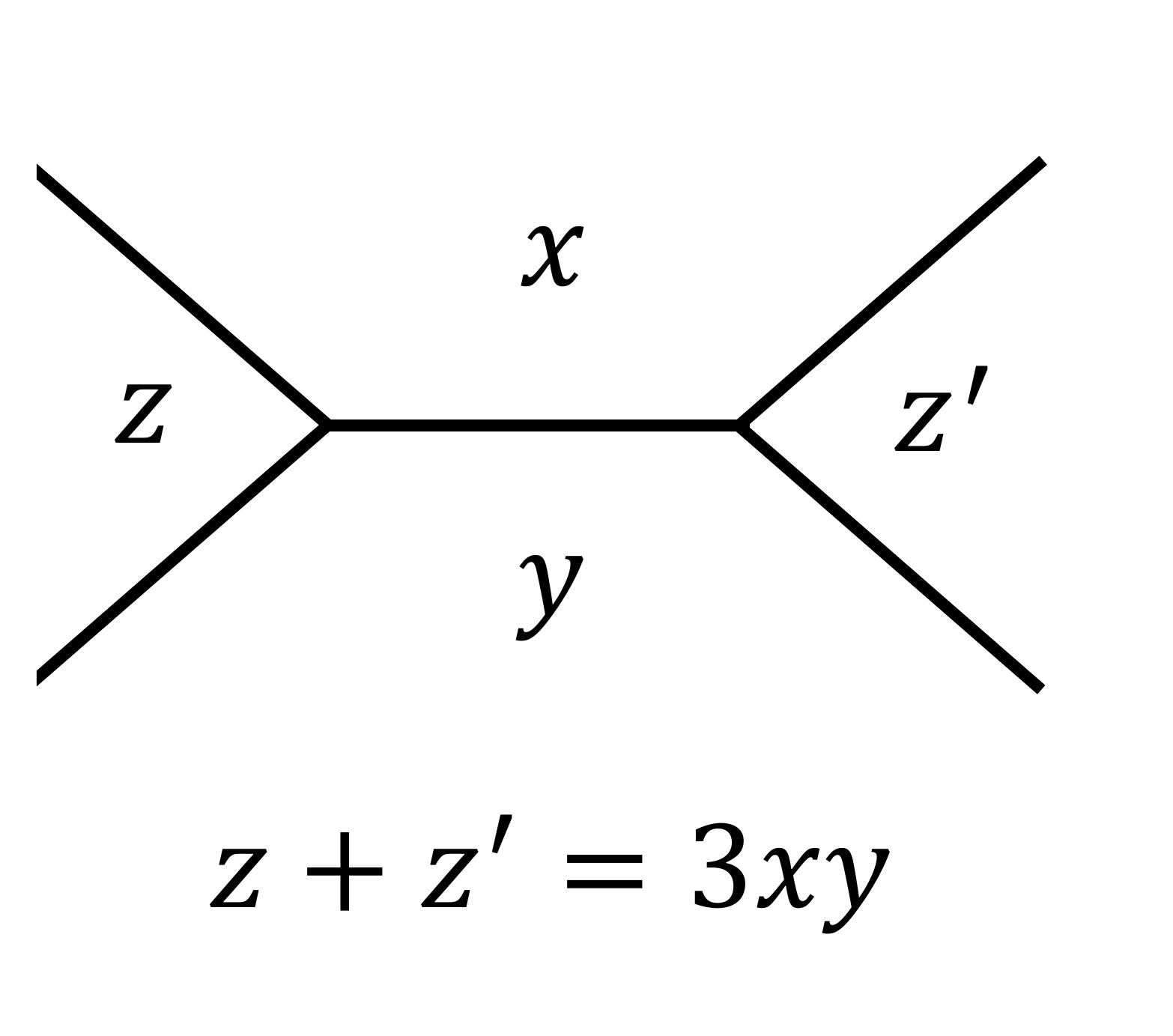} \quad \quad  \includegraphics[width=52mm]{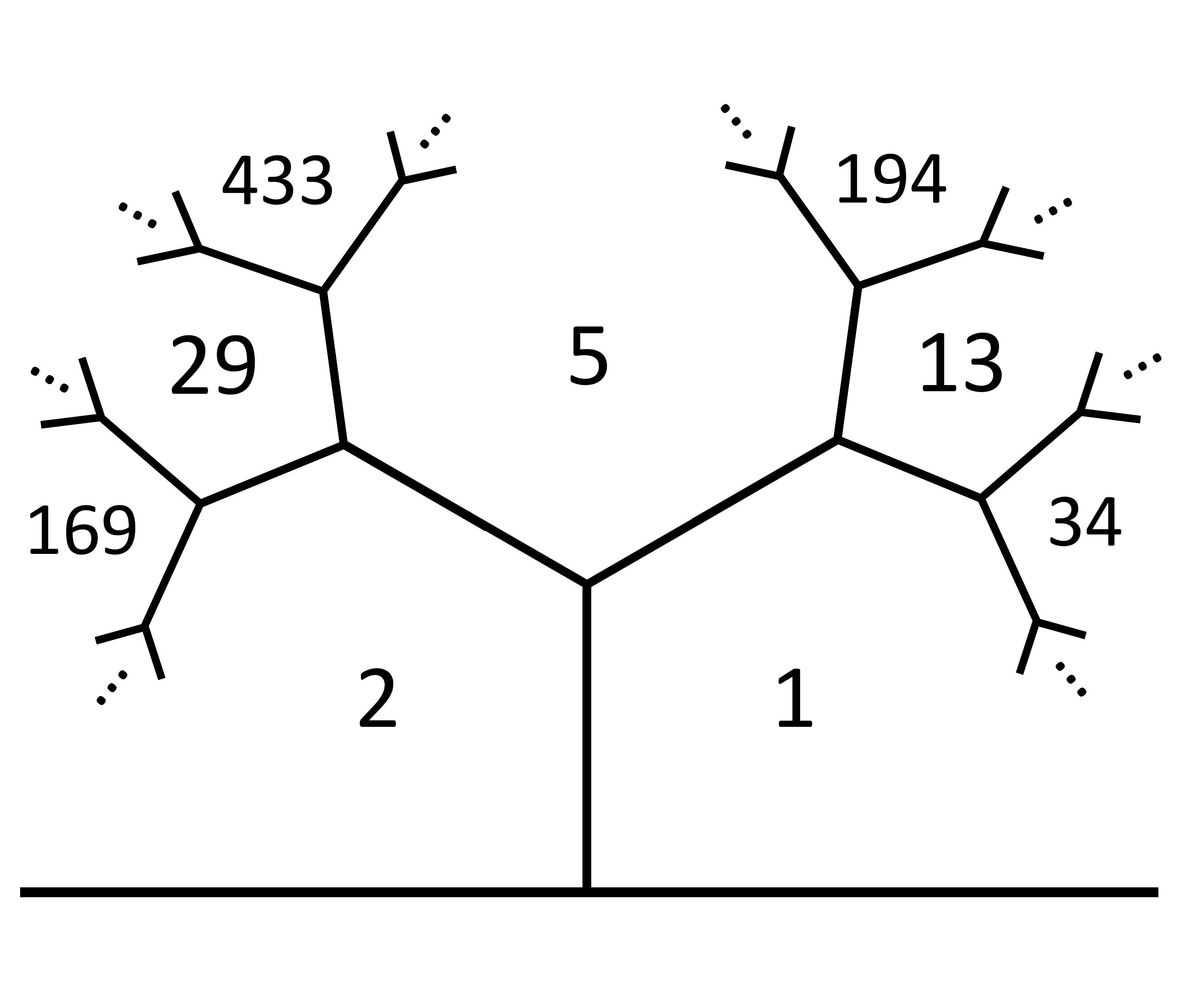} \end{center}
 \caption{\small Vieta involution and Conway topograph of Markov triples}
\end{figure}

\section*{Shadow Markov and Mordell triples}

Let $$\mathbb D(\mathbb Z)=\{a+b\varepsilon, \,\, a,b \in \mathbb Z, \, \varepsilon^2=0\}$$ be the commutative ring of integer dual numbers. By analogy with the Gaussian integers $a+bi, \, i^2=-1, \, a,b \in \mathbb Z$ we will call them {\it Clifford integers.} The invertible elements (units) in $\mathbb D(\mathbb Z)$ have the form $\pm 1+ b\varepsilon$ with any $b \in \mathbb Z$.

Ovsienko \cite{O1} proposed the shadow version of Markov triples as the solutions over Clifford integers of the following version of the Markov equation  
 \begin{equation}
\label{dM}
X^2 + Y^2 + Z^2=(3-2\varepsilon)XYZ, \quad X,Y,Z \in \mathbb D(\mathbb Z)
\end{equation}
as $PSL_2(\mathbb Z)$-orbit of the initial triple of units $X=1, \, Y=Z=1+\varepsilon$ (the motivation of this choice see in \cite{O1}).

\begin{figure}[h]
	\centering
	 \includegraphics[height=65mm]{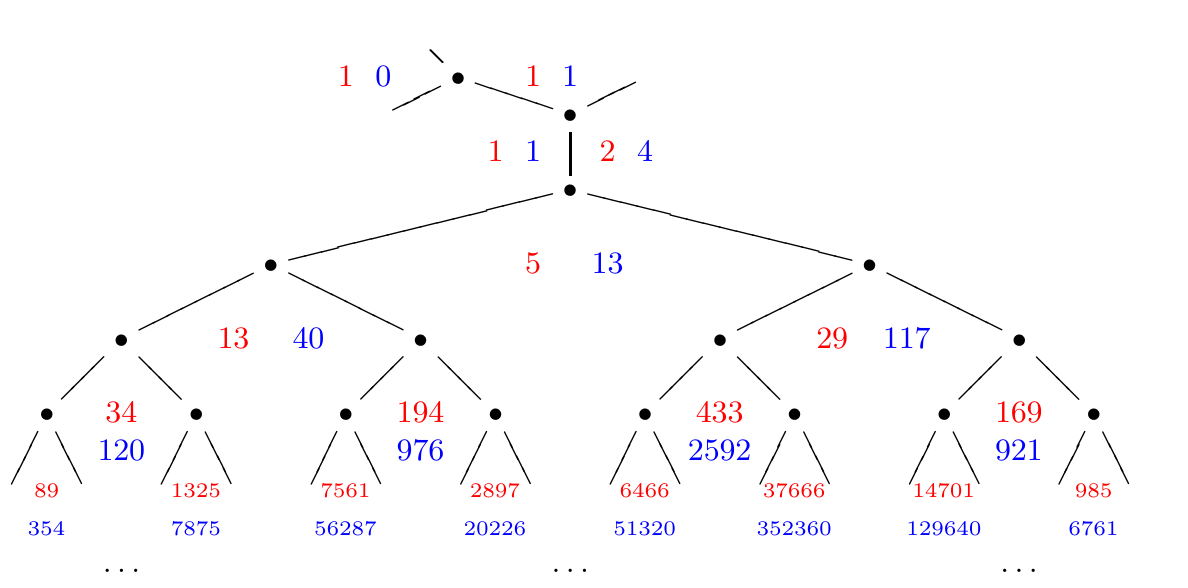}
	\caption{The shadow Markov tree (from \cite{O1})}
	\end{figure}

Ovsienko observed that the shadow companion of the very left (Fibonacci) branch of the Markov tree
$$1, 4, 13, 40, 120, 354, 1031, 2972, 8495, . . .$$
is the known sequence  A238846, which is the convolution of two bisections of the Fibonacci sequence, and asked for possible meaning of the other shadows (see Fig. 4).
This is a very interesting question, which is still largely open.

We will discuss a simpler Diophantine equation
\beq{mordell}
x^2 + y^2+ z^2 = 2xyz+1
\eeq
which was studied by Mordell \cite{Mordell} and can be viewed as a ``solvable" modification \footnote{It is very interesting to read Mordell's comments on page 505 in \cite{Mordell} about what the solvability could mean for the Diophantine equations.}  of the Markov equation. Recently it was realised that it is closely related to the two-valued formal group in complex $K$-theory studied by Buchstaber and Novikov (see \cite{BV}).

Mordell showed that this equation has integer solutions of two types: $$(I): \,\, x=1,\, y=z$$ and 
\beq{type2}
(II): \,\,2x=\xi^{a}+\eta^{a}, \, 2y=\xi^{b}+\eta^{b}, \, 2z=\xi^{c}+\eta^{c},
\eeq
where $a,b,c\in \mathbb Z$ satisfy $a+b=c$ and $\xi=p+q\sqrt{d}, \, \eta=p-q\sqrt{d},$
$p, q$ are solutions of {\it Pell's equation}
$$
 p^2-d q^2=1,
$$
where $d \in \mathbb N$ is not total square.
It is well known that all positive integer solutions of this classical equation (mistakenly ascribed to Pell by Euler) can be given by
$$
p+q\sqrt{d}=(p_0+q_0\sqrt{d})^k, \,\, k \in \mathbb N,
$$
where $(p_0, q_0)$ is the minimal positive (fundamental) solution, which can be found from the continued fraction expansion of $\sqrt{d}$ (see e.g. \cite{LeVeque}).
One can also use Conway topograph to describe the solutions of Pell's equation, as it is nicely explained by Weissman \cite{Weis}.

In particular, when $d=2$ the fundamental solution of Pell's equation $p^2-2q^2=1$ is $p_0=3,\, q_0=2$, which leads to the Mordell triples given by (\ref{type2}) with 
$
\xi=3+2\sqrt{2}, \,\,\, \eta =3-2\sqrt{2}.
$

As in the Markov case, we have the action of $PSL_2(\mathbb Z)$ generated by the Vieta involution 
$$
(x,y,z)\rightarrow (x,y, 2xy-z)
$$
and by cyclic permutations of the variables. However, unlike in the Markov case, we see that $PSL_2(\mathbb Z)$ has infinitely many orbits here.
In particular, we have a special type I orbit consisting of the single point $(1,1,1)$ fixed under the action of $PSL_2(\mathbb Z)$.

Let us look at the shadows $X=1+a\varepsilon, \, Y=1+b\varepsilon, \, Z=1+c\varepsilon$ of this special orbit.
The corresponding orbit on the Conway topograph satisfies the Vieta rule $Z+Z'=2XY$, giving
$$
c+c'=2(a+b),
$$
which is none other than Conway's arithmetic progression rule! Thus we have the following nice surprise.

\begin{Theorem}
Shadow $PSL_2(\mathbb Z)$-orbits of the special Mordell triple $(1,1,1)$ can be naturally visualized by Conway topographs of the values of the binary quadratic forms.
\end{Theorem}

%

\begin{figure}[ht]
	\centering
	 \includegraphics[height=40mm]{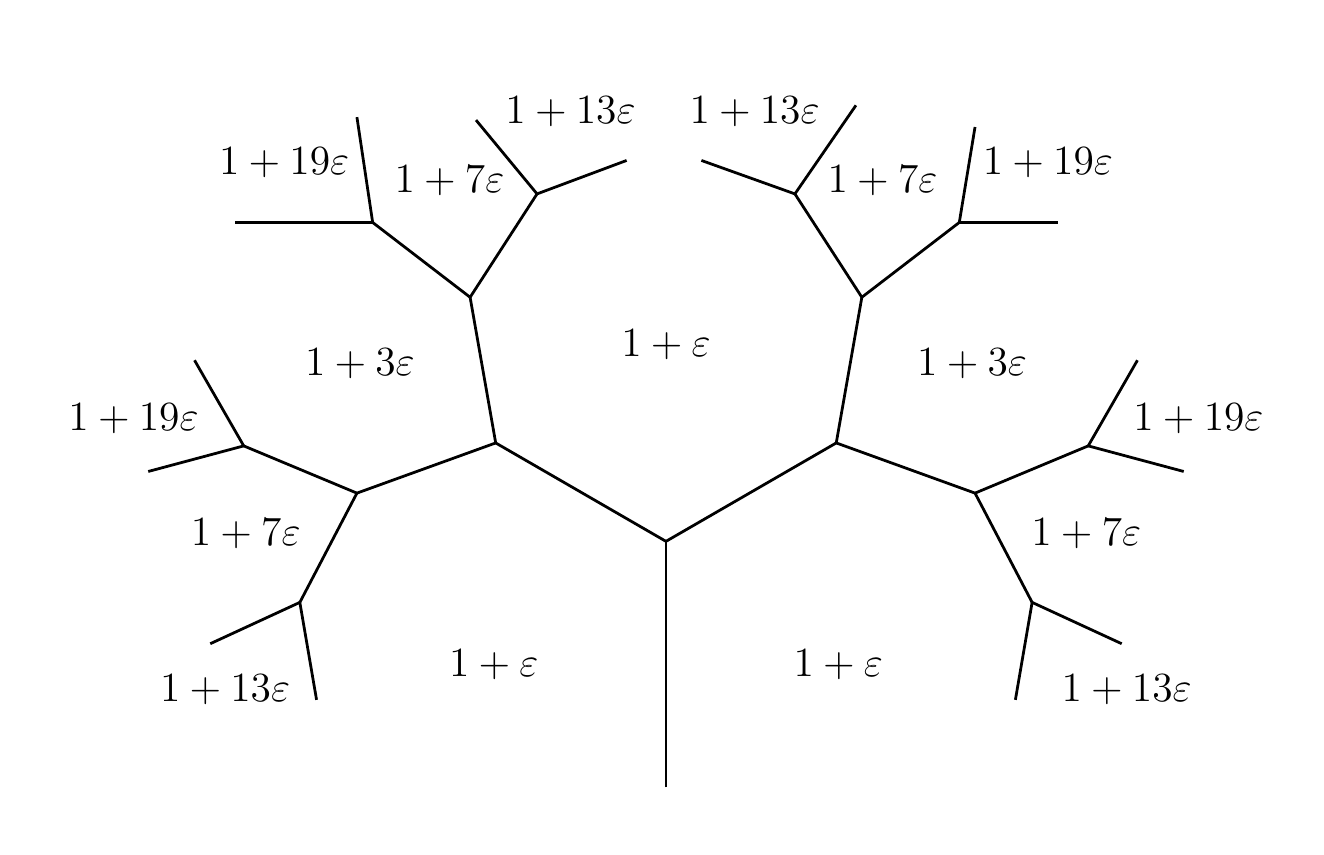} \includegraphics[height=40mm]{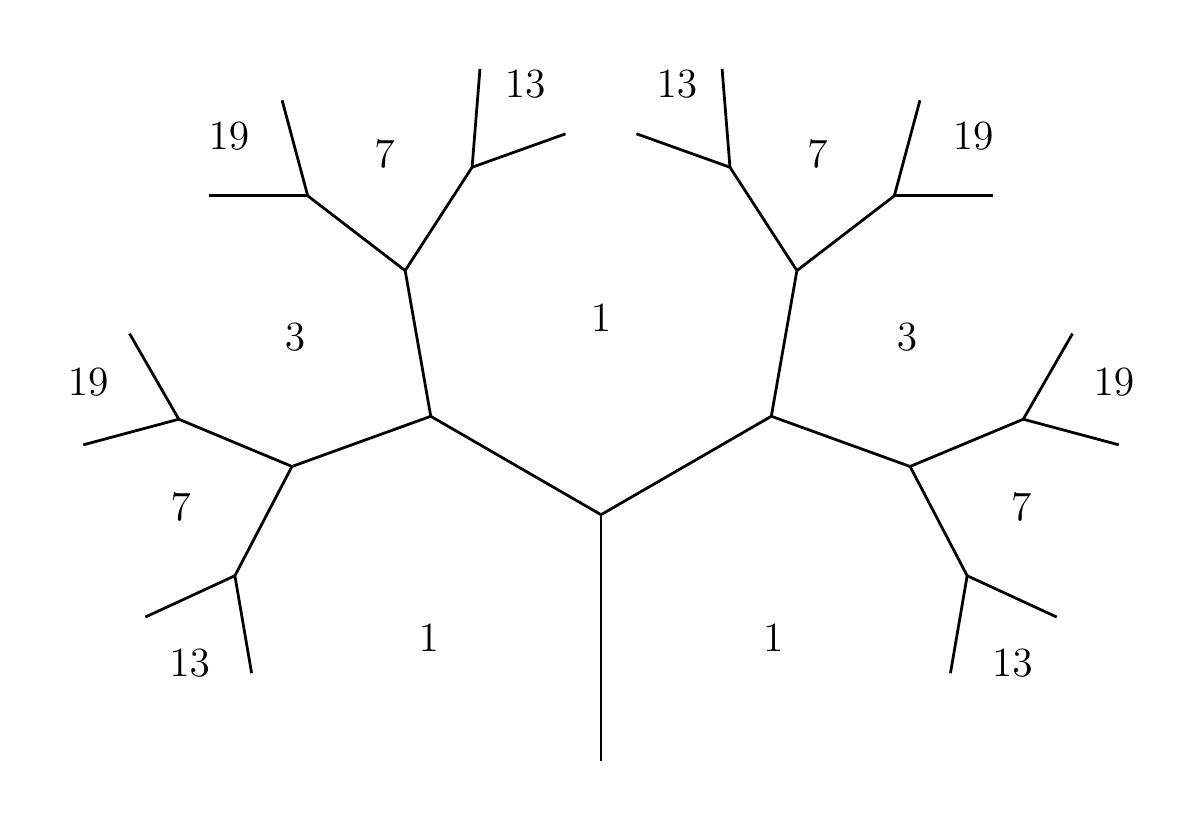}
	\caption{Shadow of special orbit $(1,1,1)$ and Conway topograph of $Q=x^2+xy+y^2$}
	\end{figure}

Note that this is true for any coefficients of the binary quadratic forms and any Mordell triples, not necessarily integer.


Before considering the shadows of other Mordell triples, let us look at the underlying geometry.
Geometrically Mordell equation (\ref{mordell}) determines a particular affine realisation of the classical {\it Cayley's nodal cubic surface} with the maximal number (which is 4) of conical singularities at the points $(1,1,1), (1, -1, -1),(-1,1,-1), (-1,-1,1)$ (see the real version on Fig 6). 

\begin{figure}[ht]
\begin{center}
\includegraphics[width=45mm]{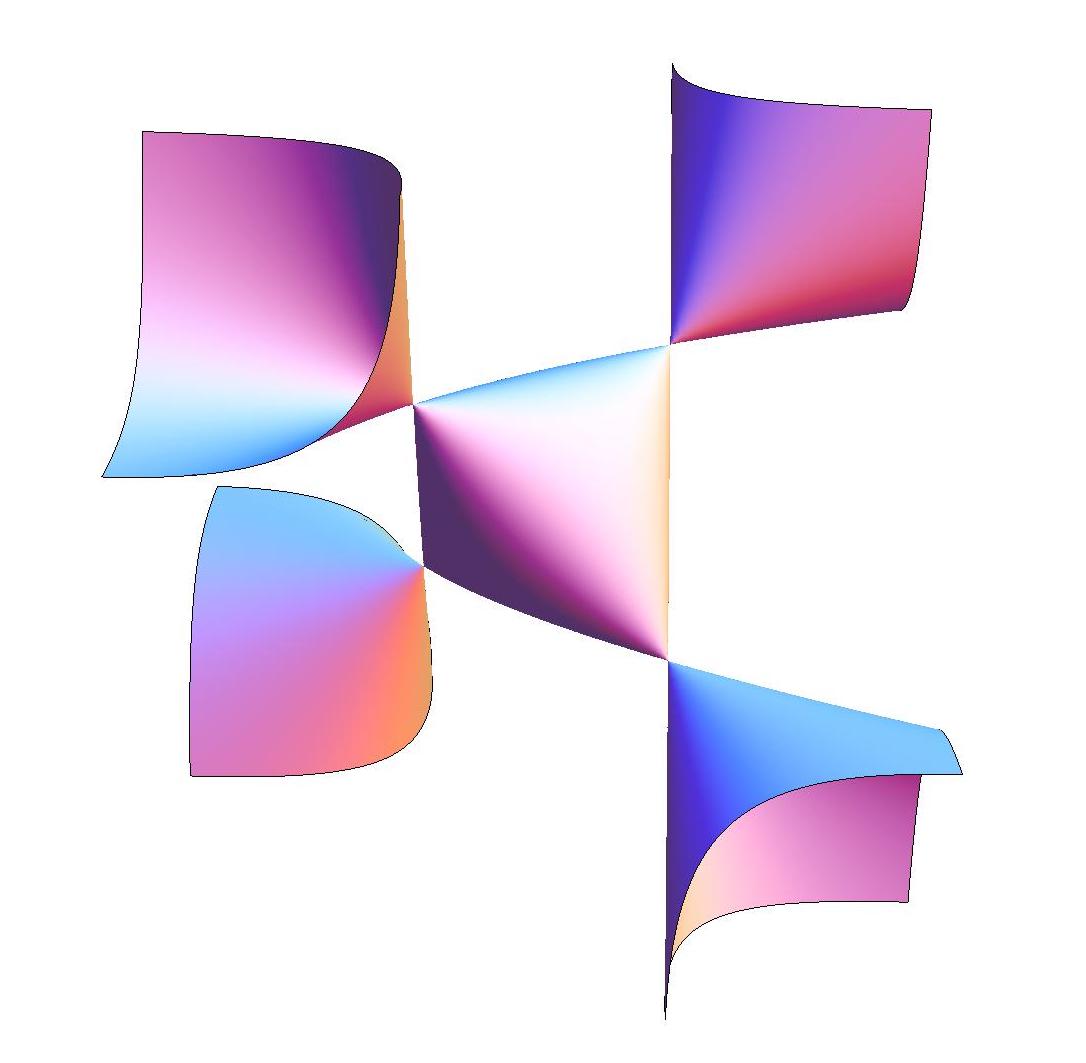}
\caption{Real Cayley-Mordell surface}
\end{center}
\end{figure}

The following observation, going back to Mordell, allows to linearise the Mordell equation. 

\begin{Theorem} \cite{Mordell}
The real Mordell triples $(x,y,z)$ with $1\leq x,\, y \leq z$ can be parametrised as
\beq{param}
x=\cosh u, \, y=\cosh v, \, z=\cosh w
\eeq
with real $u, v, w \geq 0$ satisfying the Euclid relation
\beq{lin}
w=u+v.
\eeq
\end{Theorem}


Now we can use the fact that any analytic relation can be extended to the dual numbers using the formula
$$
f(a+b\varepsilon)=f(a)+bf'(a)\varepsilon
$$
since $\varepsilon^2=0.$
This leads us to the following shadow Mordell triples:
\beq{123}
X=\cosh u+\alpha \varepsilon \sinh u, \,\,
Y=\cosh v+\beta \varepsilon \sinh v, \,\,
Z=\cosh w+\gamma \varepsilon \sinh w,
\eeq
where parameters satisfy the Euclid relation 
\beq{4}
w=u+v,\quad \gamma=\alpha+\beta.
\eeq

\begin{Theorem}
The formulae (\ref{123}),(\ref{4}) describe all solutions of the Mordell equation
\beq{shmoreq}
X^2+Y^2+Z^2=2XYZ+1
 \eeq
over real dual numbers having real parts $1\leq x,\, y \leq z$ with the exception of the special triple $x=y=z=1$, for which the shadows were described above.
\end{Theorem}

Firstly note that $X=x+\tilde x \varepsilon, \, Y=y+\tilde y \varepsilon, \, Z=z+\tilde z \varepsilon$ satisfy (\ref{shmoreq}) iff
\beq{system}
\begin{cases}
x^2+y^2+z^2=2xyz+1\\
(x-yz)\tilde x+(y-xz)\tilde y + (z-xy)\tilde z=0.
\end{cases}
\eeq
Let $X=\cosh u + \tilde x\varepsilon, \, Y=\cosh v + \tilde y\varepsilon, Z=\cosh w + \tilde z\varepsilon$ be a shadow of triple (\ref{param}) with $w=u+v.$
Then using the addition formula for $\cosh$-function 
we can rewrite the second equation in (\ref{system}) as
$$
\sinh v \sinh w \, \tilde x +\sinh u \sinh w \, \tilde y-\sinh u \sinh v \, \tilde z=0.
$$
Dividing this equation by $\sinh u \sinh v\sinh w$, assumed to be non-zero, we see that the parameters 
$\alpha=\tilde x/\sinh u, \, \beta=\tilde y/\sinh v, \, \gamma=\tilde z/\sinh w$
satisfy the relation $\gamma=\alpha+\beta$. This means that if none of $x,y,z$ is 1 then the corresponding shadows indeed have the form  (\ref{123}),(\ref{4}).
For the Mordell triples $(1,y,y)$ with $y\neq 1$ the shadows have the form $(1,y+\tilde y \varepsilon, y+\tilde z \varepsilon)$ with arbitrary $\tilde y, \tilde z,$ which agrees with (\ref{123}), (\ref{4}) since $\alpha$ could be arbitrary because $\sinh u=0.$ The case of the special triple $(1,1,1)$ is the only exceptional one here, but it was already discussed before.

Thus, the shadows of Mordell's triples (\ref{123}) on the Conway topograph are determined by the pairs of Euclid's triples (\ref{4}).
In particular, we have a natural choice when these triples are proportional: $$\alpha=\mu u, \, \beta=\mu v, \, \gamma=\mu w, \,\,\, u\pm v\pm w=0,$$ leading to the following shadows, which we will call {\it principal}:
$$
X=\cosh (u+\alpha \varepsilon)=\cosh u+\mu u \varepsilon \sinh u,
$$
$$
Y=\cosh (v+\beta\varepsilon)=\cosh v+ \mu v \varepsilon \sinh v, 
$$
$$
Z=\cosh (w+\gamma \varepsilon)=\cosh w+ \mu w \varepsilon \sinh w.
$$

Let us come back now to the integer Mordell triples of type II
\beq{mor}
x=\frac{1}{2}(\xi^{a}+\eta^{a}), \, y=\frac{1}{2}(\xi^{b}+\eta^{b}), \, z=\frac{1}{2}(\xi^{c}+\eta^{c}),
\eeq
where integers $a,b,c$ satisfy $a+b=c$ and $\xi=p+q\sqrt{d}, \, \eta=p-q\sqrt{d}$
with $p,q$ being solution of Pell's equation
$p^2-d q^2=1.$


\medskip

\begin{Theorem} The following formulae provide an explicit form of the principal integer shadows 
of Mordell triples (\ref{mor}):
\beq{shmor}
\tilde x=\frac{1}{2\sqrt{d}}ma(\xi^{a}-\eta^{a}),\,\,\,
\tilde y=\frac{1}{2\sqrt{d}}mb(\xi^{b}-\eta^{b}),\,\,\,
\tilde z=\frac{1}{2\sqrt{d}}mc(\xi^{c}-\eta^{c}),
\eeq
where  $m, a,b,c  \in \mathbb Z, \, a+b=c.$
\end{Theorem}

Indeed, since $\xi\eta=(p+q\sqrt{d})(p-q\sqrt{d})=p^2-dq^2=1$ then the equality $\cosh u=\frac{1}{2}(\xi^{a}+\eta^{a})$ implies that $u=\pm a\log\xi$ and $\sinh u=\pm \frac{1}{2}(\xi^{a}-\eta^{a})$. Now we need only to choose the values of the parameter $\mu$ in the formula for principal shadows to make them integer, which leads to formulae (\ref{shmor}).

In particular, for Pell's equation $p^2-2q^2=1$ with $p=3, q=2, m=1$ we have the following Conway topographs of the corresponding Mordell triples and their shadows:

\begin{figure}[htbp]
\begin{center}
\includegraphics[width=62mm]{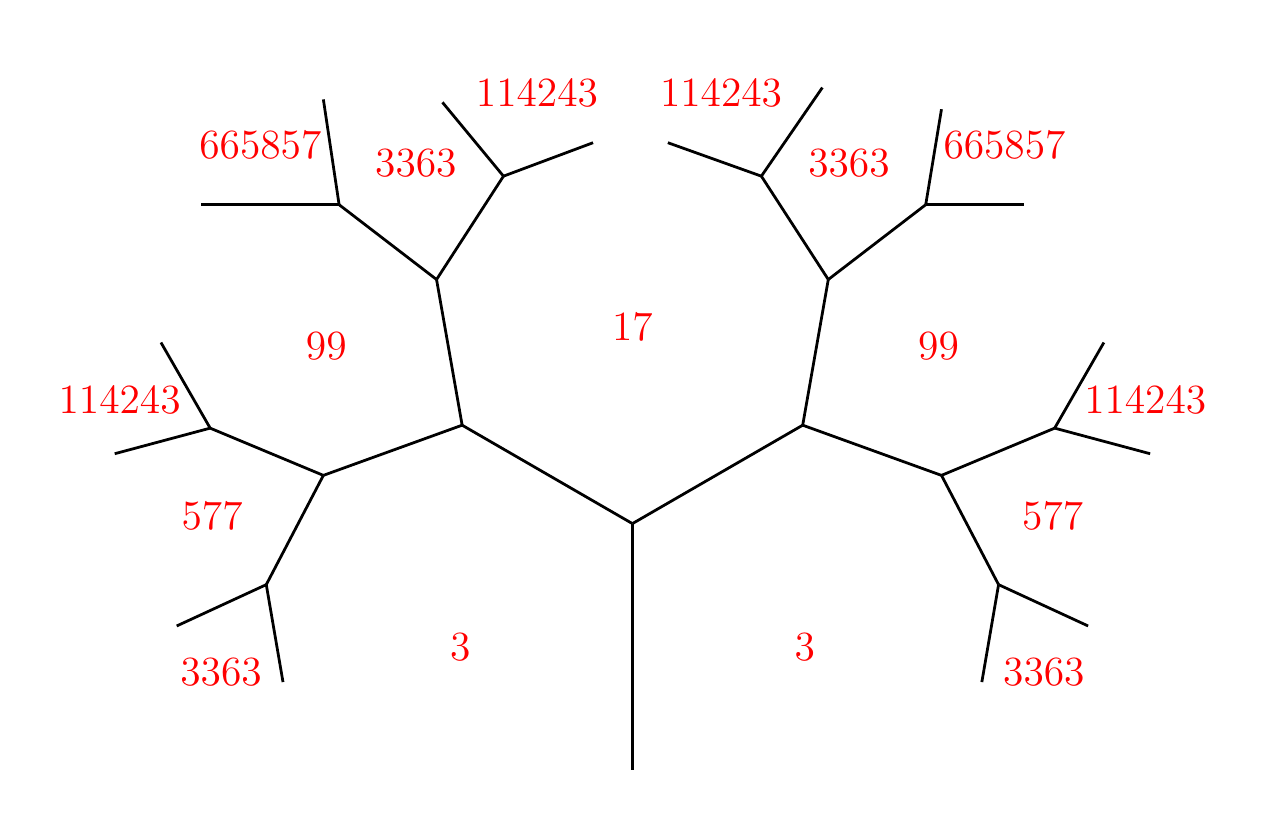} \includegraphics[width=62mm]{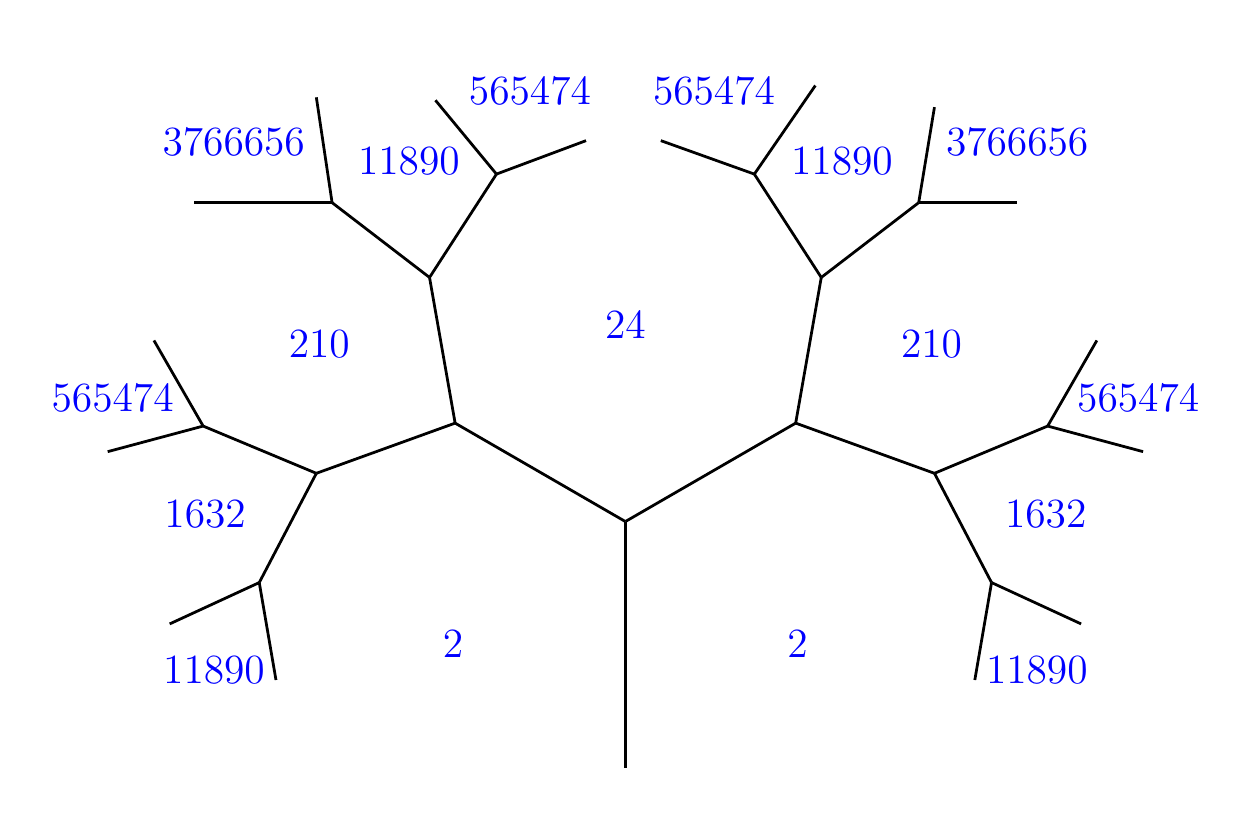}
\caption{Mordell tree and its shadow}
\end{center}
\end{figure}

Note that in contrast to Markov equation, where we have only the local Vieta rule on the topograph to generate the solutions, in Mordell case we have also an explicit formula (\ref{shmor}) for the Mordell triples and their shadows in terms of the corresponding Euclid triples (cf. Mordell's comments in \cite{Mordell}).

\section*{Shadow growth on Conway topograph}

As we can see in Figure 7, on the Conway topograph the shadows grow faster than Mordell triples. Let us look at their growth in more detail.

Consider the corresponding {\it Euclid tree} of triples $a,b,c \in \mathbb Z$ with $a+b=c,$ parametrising Mordell triples (\ref{mor}).
Note that $2x(a):=\xi^{a}+\eta^{a}=2x(-a)$, so we can restrict ourselves with positive $a,b,c$ (see Fig. 8). Going down this tree can be viewed as an application of the classical Euclid algorithm of finding the greatest common divisor of $a$ and $b$ (assumed here to be 1).

\begin{figure}[htbp]
\begin{center}
\includegraphics[width=90mm]{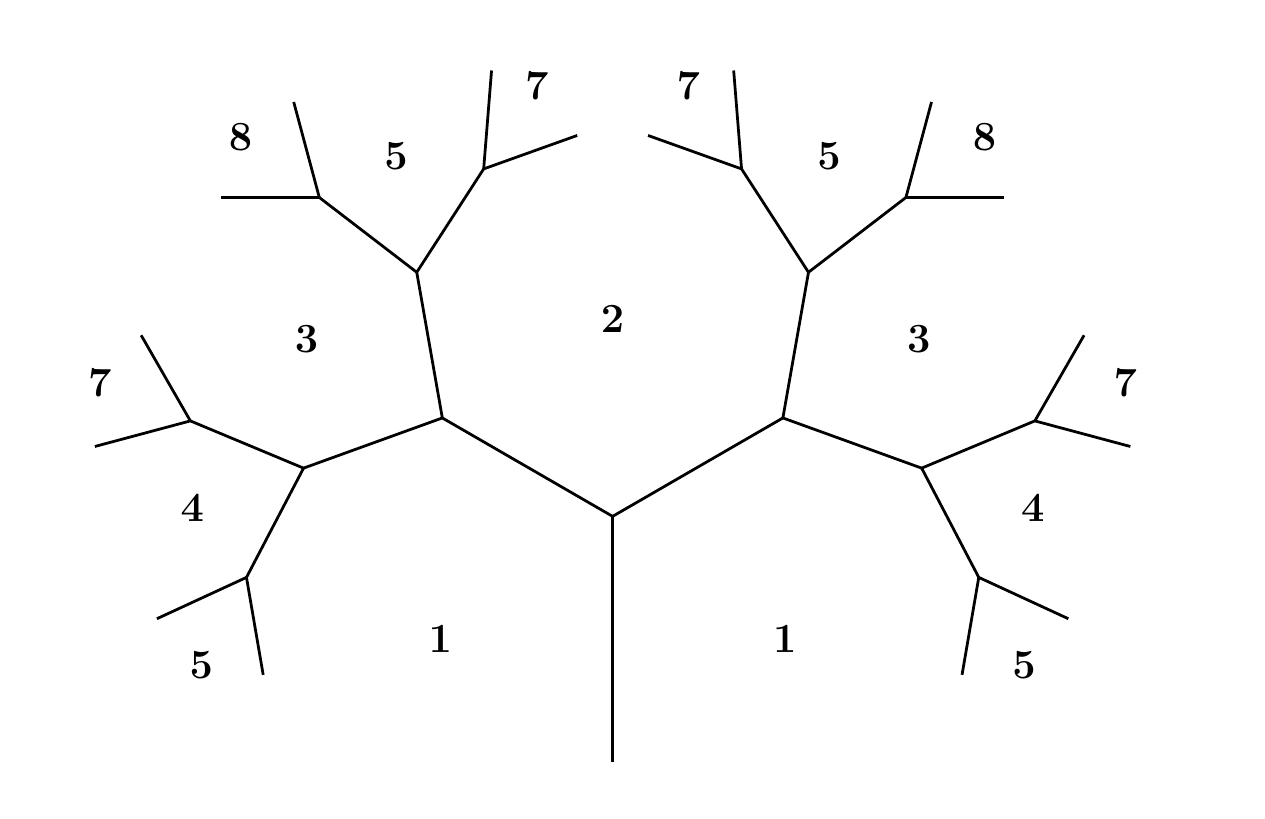}
\caption{Euclid tree}
\end{center}
\end{figure}

\begin{Theorem} 
On the Conway topograph the growth of Mordell numbers and their principal shadows (\ref{shmor}) are related asymptotically for large $a$ by
$$
\tilde x(a)\sim C\,a\, x(a), \quad C=m/\sqrt{d},
$$
where $a$ is the corresponding entry on the Euclid tree.

In other words, asymptotically we have the shadow growth
$$
\tilde x \sim c x \ln x, \quad c=m/\sqrt{d}\ln \xi.
$$
\end{Theorem}

Indeed, we know that $\xi=p+q\sqrt{d}>1, \, \eta=p-q\sqrt{d}<1$ since $\xi\eta=p^2-dq^2=1.$ This implies that for large positive $a$ we asymptotically have
$2x(a)=\xi^a+\eta^a \sim \xi^a$ and thus on the Conway topograph the shadows grow as
$$
2\tilde x(a)=ma(\xi^a-\eta^a)/\sqrt{d} \sim ma \, \xi^a/\sqrt{d}\sim C\,a\, x(a).
$$
Here we used the standard notation $f(a)\sim g(a)$ when $\lim_{a\to\infty}\frac{f(a)}{g(a)}=1.$

Let us look now more closely at the growth along the paths on the Conway topograph.
Using the Farey tree we can label the directed infinite paths by $\xi\in \mathbb RP^1$, where $\xi$ is the limit of Farey fractions along the path (see Fig. 9).

\begin{figure}[h]
\begin{center}
\includegraphics[width=59mm]{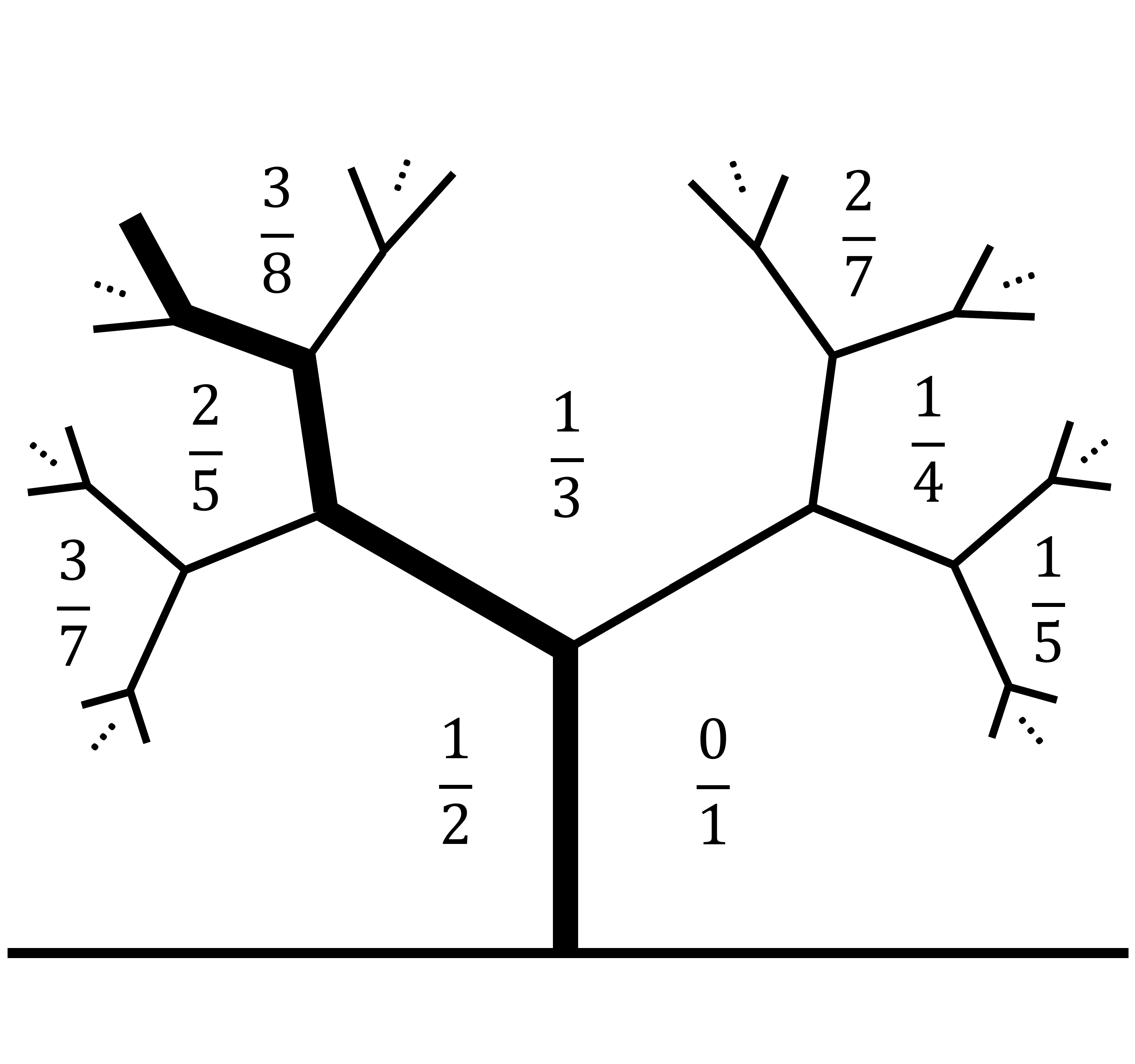}  \quad \includegraphics[height=57mm]{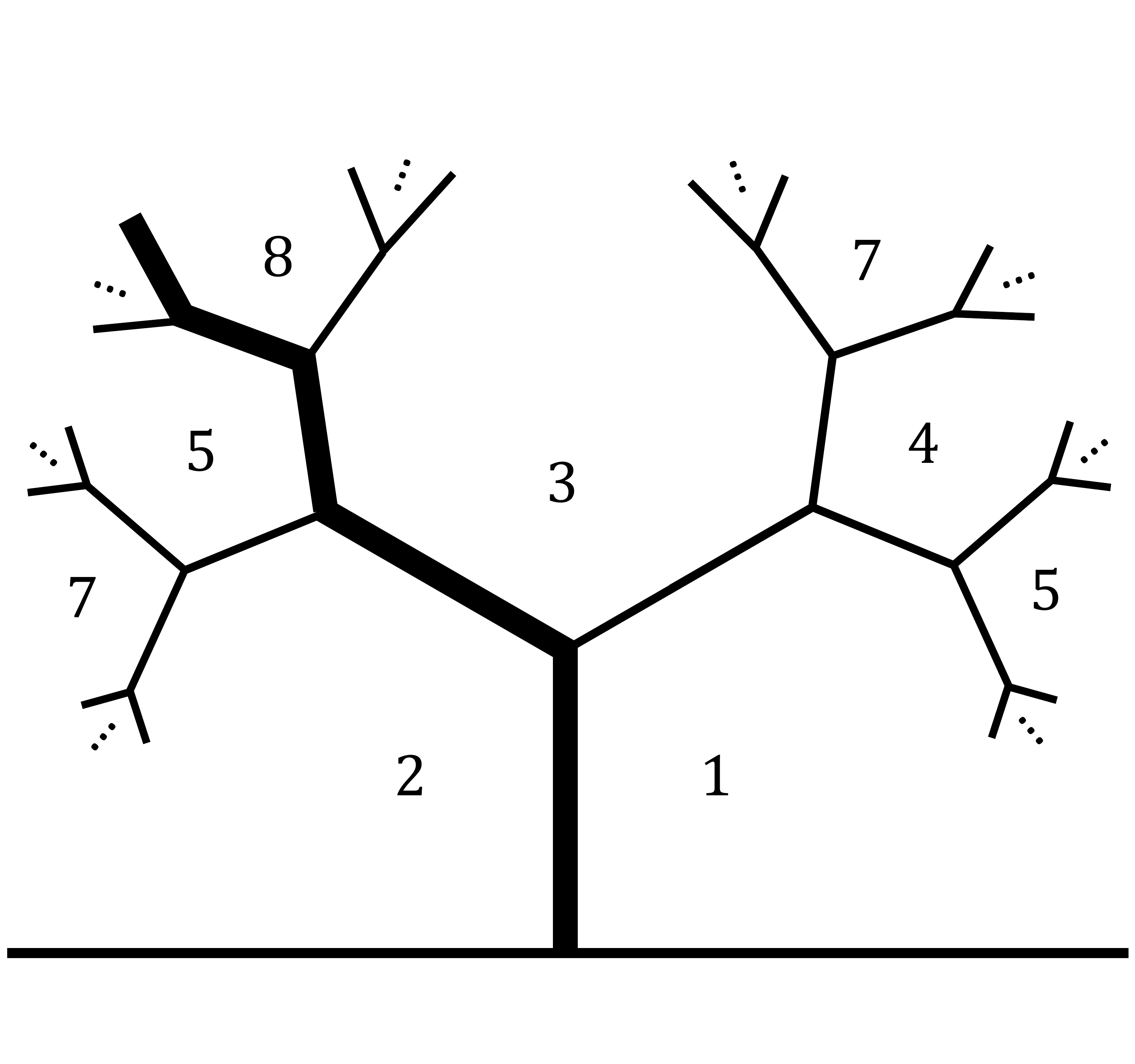}
\caption{\small Farey and Euclid trees with the ``golden" Fibonacci path}
\end{center}
\end{figure}

Following \cite{SV1} define the {\it Lyapunov function} $\Lambda(\xi)$ describing the growth of the Euclid triples $(a_n(\xi),b_n(\xi),c_n(\xi))$ along the path $\gamma(\xi):$
\beq{lam}
\Lambda(\xi):=\limsup_{n\to\infty}\frac{\ln |a_n(\xi)|}{n}
\eeq
(here $a_n(\xi)$ can be replaced by $b_n(\xi)$ or $c_n(\xi)$).

\bigskip

\begin{Theorem} \cite{SV1} The function $\Lambda(\xi)$ has the following properties:
\begin{itemize}
\item $\Lambda(\xi)$ is defined for all $\xi \in \mathbb{R}P^1$ and is $GL_2(\mathbb Z)$-invariant:
$$
\Lambda\left(\frac{a\xi+b}{c\xi+d}\right)=\Lambda(\xi), \quad \xi \in \mathbb{R}P^1,\, \begin{pmatrix}
  a & b \\
  c & d \\
\end{pmatrix} \in GL_2(\mathbb Z)
$$

\item
$\Lambda(\xi)$ vanishes almost everywhere, but the Hausdorff dimension of its support is 1.

\item
Lyapunov spectrum $Spec_E:=\{\Lambda(\xi), \,\, \xi \in  \mathbb{R}P^1\}$ of the Euclid tree is 
$$
Spec_E=[0, \ln \varphi],
$$
where $\varphi=\frac{1+\sqrt 5}{2}$ is the golden ratio.

\end{itemize}

\end{Theorem}


Using this function we can describe the growth of both Mordell triples and their shadows.

The special orbit $(1,1,1)$ is stationary for the modular dynamics, so it has zero growth, but, as we have seen, the shadows are the values of binary quadratic forms.
Their growth on the Conway topograph was studied by Spalding and the author in \cite{SV2}, which leads to the following result.

Consider the growth on the Conway topograph of the $PSL_2(\mathbb Z)$-orbit of $X=1+a\varepsilon, \, Y=1+b\varepsilon, \, Z=1+c\varepsilon.$
Let  $$Q(x,y)=ax^2+(c-a-b)xy+by^2$$ be the corresponding binary quadratic form and 
$$D(a,b,c)=a^2+b^2+c^2-2ab-2ac-2bc$$ be its discriminant, which we assume for simplicity to be nonzero. When $D>0$ the quadratic form is indefinite, otherwise it is either positive, or negative definite \footnote{On the corresponding real projective plane the equation $D(a,b,c)=0$ defines a conic with the disc $D<0$ being the celebrated Cayley-Klein model of the hyperbolic plane.}. As a corollary of the results of \cite{SV2} we have the following

\begin{Theorem} When $D(a,b,c)<0$ then the growth of the shadows of the special Mordell triple $(1,1,1)$ along the path $\gamma(\xi)$ satisfies
$$
\limsup_{n\to\infty}\frac{1}{n}\ln |a_n(\xi)|=2\Lambda(\xi).
$$ 
When $D(a,b,c)>0$ then the same is true except for $\xi=\alpha, \bar\alpha$, which are roots of quadratic equation $Q(\xi,1)=0,$ defining the ends of the corresponding {\it Conway river}, where the growth is zero.
\end{Theorem}

Recall that the {\it Conway river} \cite{Conway} is the infinite path separating on Conway topograph positive and negative values of the indefinite binary quadratic forms (see Fig. 10 and more details in \cite{SV2}).

 \begin{figure}[h]
\begin{center}
\includegraphics[height=75mm]{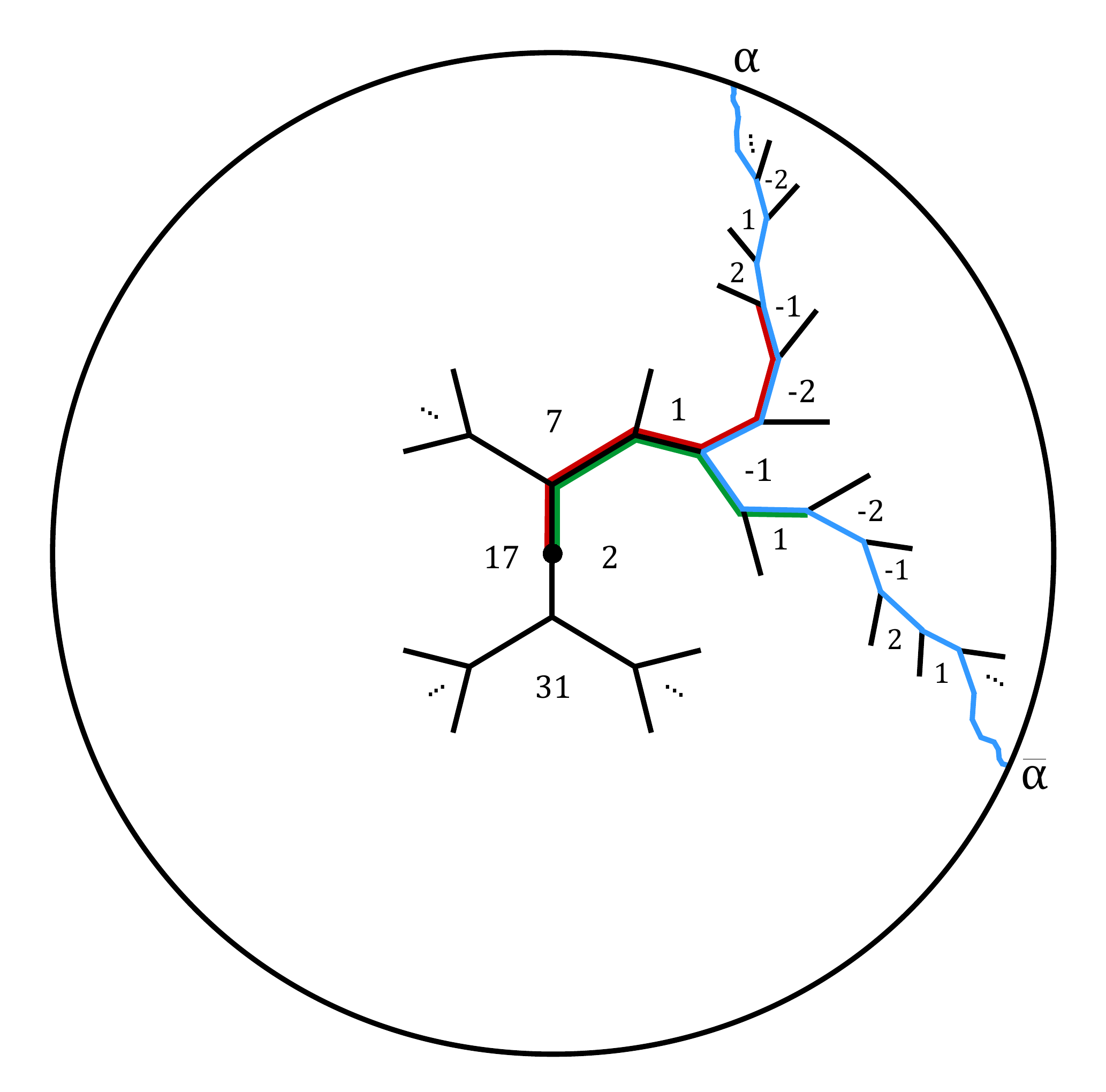} 
\caption{\small Conway river for $Q=17x^2-12xy+2y^2.$}
\end{center}
\end{figure}

For the type II Mordell triples we can describe the relative growth of the shadows along the path $\gamma(\xi)$.
Let $x_n(\xi)$ and $\tilde x_n(\xi)$ be the Mordell numbers (\ref{mor}) and their shadows (\ref{shmor}) along the path $\gamma(\xi)$.

\begin{Theorem}  The relative growth of the principal shadows $\tilde x_n(\xi)$ of the type II Mordell numbers $x_n(\xi)$ along the path $\gamma(\xi)$ can be described by
$$
\limsup_{n\to\infty}\frac{1}{n}\ln \frac{\tilde x_n(\xi)}{x_n(\xi)}=\Lambda(\xi).
$$
\end{Theorem}

Indeed, from the explicit formulae (\ref{mor}), (\ref{shmor}) it follows that $$
\limsup_{n\to\infty}\frac{1}{n}\ln \frac{\tilde x_n(\xi)}{x_n(\xi)}=
\limsup_{n\to\infty}\frac{\ln a_n(\xi)}{n}=\Lambda(\xi).
$$

Note that the condition on the shadows being principal can be removed here. For the generic shadows of the orbit $(1,1,1)$ this limit is $2\Lambda(x)$.

I believe that similar results about shadow growth hold also for Ovsienko's shadow Markov numbers, but this is still to be justified.

Another natural problem would be to study similar questions for the elliptic version of the Mordell equation
\begin{equation}\label{jac}
X^2+Y^2+Z^2=1+2XYZ-k^2(1-X^2)(1-Y^2)(1-Z^2),
\end{equation}
corresponding to the addition formula for the Jacobi elliptic function $cn(u,k)$ (see \cite{BV}).
As it is explained in \cite{BV}, these equations have significant meaning in topology: Mordell equation is closely related to the formal group law in complex $K$-theory, while its elliptic version to the elliptic cohomology. 


\section*{Acknowledgements}

I am very grateful to Valentin Ovsienko, who has initiated my interest to the shadow dynamics, and to Andy Hone for helpful and encouraging discussions.
I would like also to thank Katie Spalding and Sam Evans for helping to prepare the figures of Conway topographs.

\end{document}